# GENERALIZATION OF USUAL CAPABILITY INDICES FOR UNILATERAL TOLERANCES

DANIEL GRAU

Laboratory of Applied Mathematics, CNRS UMR 5142
IUT de Bayonne, Université de Pau et des Pays de l'Adour,
3 Av. Jean Darrigrand, 64100 Bayonne, France
Email: daniel.grau@iutbayonne.univ-pau.fr

**Abstract :**
Capability indices are dimensionless quantities measuring the aptitude of a process to manufacture items whose characteristics must be within a specified tolerance ranges. The usual indices $C_p$, $C_{pk}$, $C_{pm}$, and $C_{pmk}$ are used for a process of normal distribution and a target located at the center of the tolerance interval. Various indices derived from the previous family allow to consider more complex situations when asymmetrical tolerances and non-normal distributions are taken into account. In this paper we study the case where a single tolerance is imposed because the shifts in the direction of this tolerance appear much more serious than in the opposite direction. We propose a family of four indices having interpretations and properties similar to those of the usual family.

## 1. Introduction

The quality of an industrial process is assured by the monitoring of one or more variables of interest. The process will be considered of good quality if it is able to produce items whose variable of interest is within the tolerances *L* and *U* (Lower and Upper tolerance limits) and close to a specified target value *T*. Accordingly, the capability indices used to measure the quality of the process are linked to the location and the dispersion of the supervised variable. For a target centered in its tolerance interval and a variable of normal distribution, the first work of Kane [9] was completed to lead to a family of four indices generally used. X being the supervised characteristic of normal distribution $N(\mu, \sigma)$, one defines

$C_p = \dfrac{U - L}{6\sigma} = \dfrac{d}{3\sigma}$, where $d = (U - L)/2$ is the half-length of the tolerance interval,

$C_{pk} = \dfrac{\min(U - \mu, \mu - L)}{3\sigma}$,

$C_{pm} = \dfrac{U - L}{6\sqrt{\sigma^2 + (\mu - T)^2}}$,

$C_{pmk} = \dfrac{\min(U - \mu, \mu - L)}{3\sqrt{\sigma^2 + (\mu - T)^2}}$.

When asymmetrical tolerances and normal data are concerned, many proposals have been made, the most coherent of which as the generalization of the usual family, is given by the family $C_p^{"}$, $C_{pk}^{"}$, $C_{pm}^{"}$, and $C_{pmk}^{"}$ introduced by Pearn, Chen, and Lin in various articles [3, 4, 12, 16, 17].

In the case of non-normal distributions, Clements [5] proposes to generalize the indices $C_p$ and $C_{pk}$ replacing the mean $\mu$ by the median $M$ and the dispersion $6\sigma$ by $Up - Lp$ where $Up$ and $Lp$ are the 99.865 and 0.135 percentiles. In the same way, Pearn and Kotz [15] provide a extension of the indices $C_{pm}$ and $C_{pmk}$. However let us note that considering the difficulty to obtain reliable estimations of $Up$ and $Lp$, many other proposals have been provided the bibliography of which can be found in the articles of Tang and Than [21], Mac Cormack, Harris, Hurwitz, and Spagon [11], or Ding [6].

The literature related to the processes limited by two tolerance limits is important, but it is not the same in the less frequent situation where only one limit is imposed. In this case two situations can be considered. The first one occurs when the variable of interest, because of its nature, cannot exceed a certain level which represents the target of the process. It is the case for example for concentricity or circularity where the observed measure is obviously positive or null, the target being equal to 0. The second one occurs when a drift of the mean in a direction appears much less serious to the user, so that he is induced to define only one single tolerance. It is this second situation which is the subject of this paper. We recall the rare indices which one finds in the literature, then propose a family of four indices having interpretations and properties, similar to the ones of the usual family $C_p$, $C_{pk}$, $C_{pm}$, and $C_{pmk}$. For normal data, the densities and the moments of the natural estimators are given. The expressions used in the case of a normal distribution can be easily generalized to the case of non-normal distributions by replacing the mean $\mu$ by the median $M$, and the natural variation $3\sigma$ by $Up - M$ or $M - Lp$. The estimation of the extreme percentiles $Up$ and $Lp$ is however far from reliable on samples of reasonable size when one uses the traditional moment estimation, as Clements does [5]. The work of Shore [19,20] enables us to obtain estimations much more reliable which will be developed in an example.

## 2. Existing indices

Starting from the definition of $C_{pk} = \dfrac{\min \left( U - \mu, \mu - L \right)}{3\sigma}$, Kane [8] defines $CPU = (U - \mu)/3\sigma$ and $CPL = (\mu - L)/3\sigma$ in order to measure the capability of a process in the unilateral tolerance situation. Let us note that the $CPU$ and $CPL$ indices do not take into account the existence of a target which, for the index $C_{pk}$, is supposedly, implicitely located at the center of the tolerance interval. When the target is not centered, Kane [8] suggests an index referred to as $C_{pk}^*$ from which he gets the indices $CPU^* = \left( U - T - \left| T - \mu \right| \right)/3\sigma$ and $CPL^* = \left( T - L - \left| T - \mu \right| \right)/3\sigma$, in the case of a unilateral tolerance. In the same way as with the usual indices, the value of these indices is equal to 0 if the previous calculus leads to a negative value. Chan, Cheng and Spring [2], in the case where one-sided tolerance is required, have suggested generalizing $C_{pm}$ to $C_{pmu}^* = \dfrac{U - T}{3\sqrt{\sigma^2 + \left( \mu - T \right)^2}}$ and

$$C_{pml}^* = \frac{T - L}{3\sqrt{\sigma^2 + \left( \mu - T \right)^2}} \, .$$

On the same principle, Vänmann [23] suggests generalizing $C_{pmk}$ to

$$C_{pmku} = \frac{U - \mu}{3\sqrt{\sigma^2 + (\mu - T)^2}} \text{ and } C_{pmkl} = \frac{\mu - L}{3\sqrt{\sigma^2 + (\mu - T)^2}}.$$ Moreover in order to generalize all

the previous suggestions, Vänmann [23] puts forward two families of indices,

$$C_{pau}(u,v) = \frac{U - \mu - u|\mu - T|}{3\sqrt{\sigma^2 + v(\mu - T)^2}} \text{ and } C_{pal}(u,v) = \frac{\mu - L - u|\mu - T|}{3\sqrt{\sigma^2 + v(\mu - T)^2}} \text{ on the one hand,}$$

$$C_{pvu}(u,v) = \frac{U - T - u|\mu - T|}{3\sqrt{\sigma^2 + v(\mu - T)^2}} \text{ and } C_{pvl}(u,v) = \frac{T - L - u|\mu - T|}{3\sqrt{\sigma^2 + v(\mu - T)^2}} \text{ on the other hand. We}$$

have $C_{pau}(0,0) = CPU$, $C_{pal}(0,0) = CPL$, $C_{pau}(0,1) = C_{pmku}$, $C_{pal}(0,1) = C_{pmkl}$,

$C_{pvu}(1,0) = CPU^*$, $C_{pvl}(1,0) = CPL^*$, $C_{pvu}(0,1) = C_{pmu}^*$, and $C_{pvl}(0,1) = C_{pml}^*$. As Vänmann

points out, the fact that there is an unilateral tolerance can be interpreted in such a way that a

shift of $\mu$ away from $T$ towards that tolerance is more serious than a shift of $\mu$ towards the

opposite side. The indices derived from $C_{pvu}(u,v)$ and $C_{pvl}(u,v)$, which are symmetrical

around the target are thus of little interest. From the properties of the estimators of the indices

$C_{pau}(u,v)$ and $C_{pal}(u,v)$, Vänmann [23] suggests using $C_{pau}(0,4)$ and $C_{pal}(0,4)$, although

these indices are not maximum when $\mu$ is on the target $T$ (fig 1). This drawback is not deemed

too serious by the author, since a shift of $\mu$ away from $T$ to the left (for $C_{pau}(u,v)$) is less

important considering the expected percentage of nonconforming than a shift of $\mu$ away from

$T$ towards $U$.

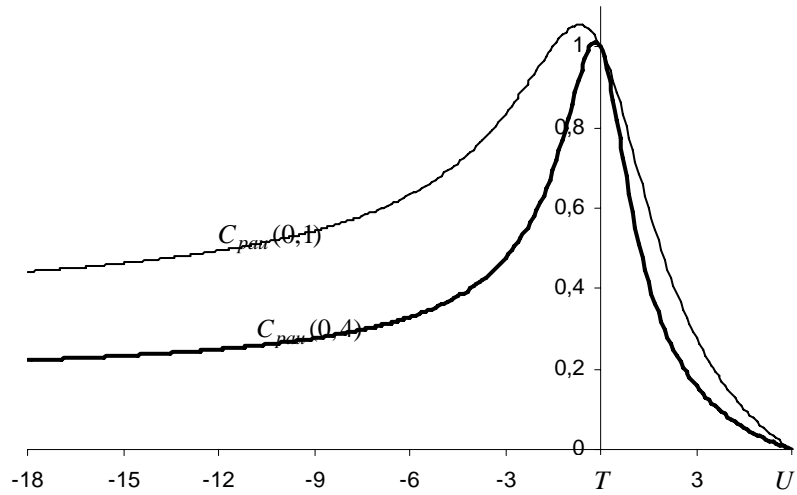

Figure 1. $C_{pau}(0,4)$ and $C_{pau}(0,1)$ as a function of $\mu$ for $U = 6$, $T = 0$, and $\sigma = 2$.

## 3. Suggestions of indices for a normal distribution

Our objective is to build four indices generalizing the usual properties and interpretations of the indices $C_p$, $C_{pk}$, $C_{pm}$, and $C_{pmk}$. Let us mention that for a normal distribution and a target centered in the tolerance interval, $C_p$ is linked to the proportion of non conforming items when the mean is on the target, and is equal to 1 for a proportion of 0,27% of non conforming. Moreover it is meant as the potential capability of the process, that is to say the maximum capability which one can obtain for a given dispersion, when the mean is on the target. The modern standard of quality deems that a process should not be considered capable if $\mu$ is far away from $T$, even if $\sigma$ is small. $C_p$ which is unrelated to $T$ does not satisfy this requirement. For this reason the indices $C_{pk}$, $C_{pm}$, and $C_{pmk}$, which take into account the location of the process mean as well as the process variability, have been introduced. The 3 of them are maximum and equal to $C_p$ when the mean is on the target. The deviation of the mean is taken linearly into account by $C_{pk}$, so that it is null when the mean reaches the tolerance limits. Thus, the ratio $C_{pk}/C_p$ allows to determine the position of the mean between the target and the tolerance towards which it deviates. $C_{pm}$ takes the deviation into account in a quadratic way. Thus, it is not null at the tolerance limits, but takes the same value. $C_{pmk}$ being the combination of $C_{pk}$ and $C_{pm}$ allows to take the deviation into account in a "quadratic" way and to obtain a null index when the mean reaches one of the tolerances. Thus, in the case of unilateral tolerance, in order to keep interpretations similar to the bilateral case, we will require

a) that the potential capability indices $C_{pu}$ and $C_{pl}$ take value 1 for a 0,135% proportion of non conforming, when the mean is on the target

b) that the indices of capability taking the position of the mean related to the target into account, are maximum and equal to potential capability when the mean is on the target.

c) that the indices $C_{pku}$ and $C_{pkl}$ decrease linearly and take value 0 when the mean reaches the tolerance limit

d) that the indices $C_{pmu}$ and $C_{pml}$ decrease in a quadratic way

e) that the indices $C_{pmku}$ and $C_{pmkl}$ decrease in a "quadratic" way and take value 0 when the mean reaches the tolerance limit.

The main difficulty in the building of indices lies in the fact that we have no knowledge of the risk taken when the mean moves away from the target in the opposite direction to the tolerance limit. However, even if there is no tolerance limit, obviously, a production manager cannot accept a too large deviation in a direction even if, a priori, it does not seem too serious to him. Thus, it appears fundamental to us that he should quantify the "not too serious". Is this twice, five times, ten times less serious? We consider thereafter that it has been decided that the risk is k time less serious. The choice of the constant k (>1) being rather approximate we require a last condition:

f) when the means deviates towards the tolerance limit, the capability indices must be independent of the choice of k.

Let us consider for the moment the case of an upper tolerance U, and note

$\alpha_u = \max\left((\mu - T)/(U - T), (T - \mu)/(k(U - T))\right)$, and $\delta_u = (U - T)\alpha_u / \sigma$. The new indices suggested are defined by :

$$C_{pu} = \frac{U - T}{3\sigma},$$

$$C_{pku} = (1 - \alpha_u) C_{pu} = \frac{U - T - \max\left(\mu - T, (T - \mu)/k\right)}{3\sigma},$$

$$C_{pmu} = \left(1 + \delta_u^2\right)^{-\frac{1}{2}} C_{pu} = \frac{U - T}{3\sqrt{\sigma^2 + \left[\max\left(\mu - T, (T - \mu)/k\right)\right]^2}},$$

$$C_{pmku} = \left(1 - \alpha_u\right) C_{pmu} = \left(1 + \delta_u^2\right)^{-\frac{1}{2}} C_{pku} = \left(1 - \alpha_u\right)\left(1 + \delta_u^2\right)^{-\frac{1}{2}} C_{pu}$$

$$= \frac{U - T - \max\left(\mu - T, (T - \mu)/k\right)}{3\sqrt{\sigma^2 + \left[\max\left(\mu - T, (T - \mu)/k\right)\right]^2}}.$$

By using notations similar to Vännman's [23], we can write the general formula

$$C_{pu}(u, v) = \frac{U - T - u\max\left(\mu - T, (T - \mu)/k\right)}{3\sqrt{\sigma^2 + v\left[\max\left(\mu - T, (T - \mu)/k\right)\right]^2}} = \frac{U - T - uA_u^*}{3\sqrt{\sigma^2 + vA_u^{*2}}}, \text{ where}$$

$A_u^* = \max\left(\mu - T, (T - \mu)/k\right)$. This expression gives the four indices for the couples $(u, v) = (0,0)$, $(1,0)$, $(0,1)$ and $(1,1)$ again. Note that the letter $u$ used in subscript is an abbreviation of the word upper and is independent of the first parameter of the indices family. The indices $C_{pu}(u, v)$ are identical Vännman's $C_{pvu}(u, v)$ indices [23] when $\mu - T > 0$, but are different when $\mu - T < 0$.

It is obvious that when $\mu = T$ and $C_{pu} = 1$, $U - T = 3\sigma$, and a proportion of 0.135% of items is thus beyond $U$, that satisfies a). If $\mu = T$, then $\alpha_u = \delta_u = 0$, thus $C_{pku} = C_{pmu} = C_{pmku} = C_{pu}$, which satisfies b). If $\mu = U$, $C_{pku} = C_{pmku} = 0$. In addition, from the previous algebraic expressions, when $\mu$ moves away from the target, it is obvious that $C_{pku}$ decreases linearly and $C_{pmu}$ in a quadratic way, hence the conditions c), d), and e). For $\mu - T > 0$, $\mu$ moves away towards $U$, and $C_{pu}(u, v)$ is thus independent of k, which satisfies the condition f). Finally, let us notice some additional properties identical to those of the usual family. From the previous algebraic relations, we have obviously $C_{pu} \geq C_{pku} \geq C_{pmku}$, $C_{pu} \geq C_{pmu} \geq C_{pmku}$, and $C_{pmku} = \frac{C_{pku} C_{pmu}}{C_{pu}}$. The $C_{pku}$ suggested, linked to $C_{pu}$, gives a precise idea of the position of the mean in the $[T; U]$ interval. Indeed, if $C_{pku}/C_{pu} = h$ and $\mu > T$, then $U - \mu = h(U - T)$. Thus for $h = \frac{1}{2}$ by example, the mean is halfway between the target and the tolerance limit. For $\mu > T$, $C_{pmu} < (U - T)/3(\mu - T)$. In particular for $C_{pmu} = 1$, $(\mu - T) < (U - T)/3$. Thus a $C_{pmu}$ value of 1 and $\mu > T$, implies that the process mean $\mu$ is in the middle third of the interval $[T; U]$. For $\mu > T$, $C_{pmku} < \left[(U - T)/3(\mu - T)\right] - 1/3$. Thus a $C_{pmku}$ value of 1 and $\mu > T$, implies that the process mean $\mu$ is in the middle fourth of the interval $[T; U]$.

To visualize the properties of the four indices, figure 2 represents the evolution of $C_{pu}$, $C_{pku}$, $C_{pmu}$, and $C_{pmku}$ according to the variations of the mean $\mu$, in the case where the deviation to the left is considered to be three times less important than to the right.

In a similar way, if the single tolerance is a lower limit $L$, we obtain the general formulation

$$C_{pl}(u,v) = \frac{T - L - u \max\left(\mu - T/k, (T - \mu)\right)}{3\sqrt{\sigma^2 + v\left[\max\left(\mu - T/k, (T - \mu)\right)\right]^2}} = \frac{T - L - uA_l^*}{3\sqrt{\sigma^2 + vA_l^{*2}}}$$ . Assume that

$\alpha_l = \max\left((\mu - T)/(k(T - L)), (T - \mu)/(T - L)\right)$, and $\delta_l = (T - L)\alpha_l/\sigma$, we find similar

algebraic expressions, $C_{pkl} = (1 - \alpha_l)C_{pl}$, $C_{pml} = \left(1 + \delta_l^2\right)^{-\frac{1}{2}}C_{pl}$, and $C_{pmkl} = (1 - \alpha_l)C_{pml}$

$= \left(1 + \delta_l^2\right)^{-\frac{1}{2}}C_{pkl} = (1 - \alpha_l)\left(1 + \delta_l^2\right)^{-\frac{1}{2}}C_{pl}$, as well as the same properties stated in the case of an upper tolerance.

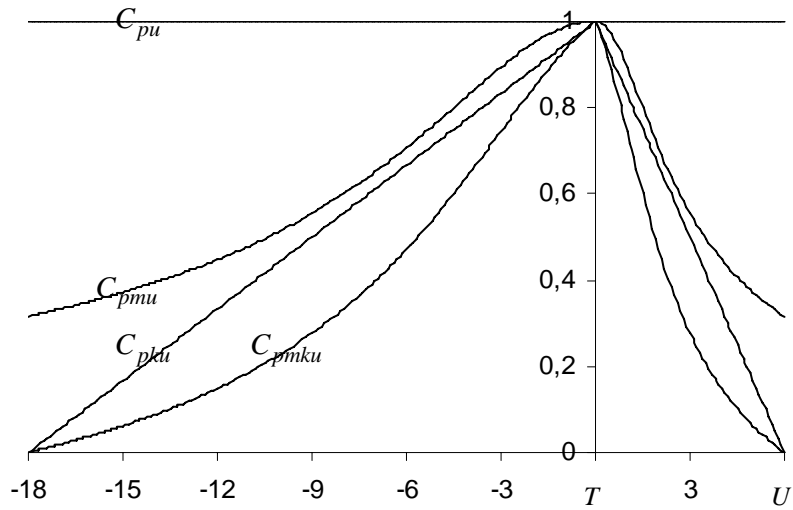

Figure 2. $C_{pu}$, $C_{pku}$, $C_{pmu}$, and $C_{pmku}$ as a function of $\mu$ for $U = 6$, $T = 0$, $\sigma = 2$, and $k = 3$.

## 4. Suggestions of indices for a non-normal distribution

As in the bilateral case, we replace the mean $\mu$ by the median $M$ and the natural variation $3\sigma$ by $U_p - M$ or $M - L_p$ according to each case. Hence the formulas

$$C_{pu}(u,v) = \frac{U - T - u \max\left(M - T, (T - M)/k\right)}{3\sqrt{\left(\frac{U_p - M}{3}\right)^2 + v\left[\max\left(M - T, (T - M)/k\right)\right]^2}},$$

and $C_{pl}(u,v) = \dfrac{T - L - u \max\left(M - T/k; (T - M)\right)}{3\sqrt{\left(\dfrac{M - L_p}{3}\right)^2 + v\left[\max\left(M - T/k; (T - M)\right)\right]^2}}$ .

When $(u,v) = (0,0)$, $(1,0)$, $(0,1)$ and $(1,1)$, the previous formulas give generalizations of the indices $C_p$, $C_{pk}$, $C_{pm}$, and $C_{pmk}$, in the case of non-normal distributions. Assuming that $\alpha_u = \max\left((M - T)/(U - T), (T - M)/k(U - T)\right)$, $\delta_u = 3(U - T)\alpha_u/(U_p - M)$, $\alpha_l = \max\left((M - T)/k(T - L), (T - M)/(T - L)\right)$, and $\delta_l = 3(T - L)\alpha_l/(M - L_p)$, we find the same algebraical expressions and thus the same properties as those stated in the case of a normal distribution. Practically, to be able to use these results we need to estimate the percentiles $M$, $U_p$, and $L_p$. The determination of a percentile is easy when the distribution from which the observations ensue is known. To identify this distribution, the most usual method consists in fitting the distribution to a member of a family covering a great number of usual distributions. To achieve this adjustment the method of the moments requires the identification of 3 or 4 parameters according to the family being used. If the estimate of the mean and of the standard deviation is rather reliable, it is no longer the same for skewness and kurtosis, which are manifestly known for their great dispersion. Thence, the resulting estimates are not at all reliable for the extreme percentiles $U_p$ and $L_p$. Thus for a non negative variable, Shore [19] suggests approaching the p percentile $Q_p$ by the relation

$$Q_p = \begin{cases} A_1\left(\dfrac{p}{1-p}\right)^{B_1} & \text{when } p < 1/2 \\[2ex] A_2 \ln\left(\dfrac{p}{1-p}\right) + B_2 & \text{when } p \geq 1/2 \end{cases} .$$

By identification of the complete moments $\mu_i = \int_0^1 x^i dF(x) = \int_\square x^i f(x) d(x)$ and partial $M_i = \int_{1/2}^1 x^i dF(x) = \int_{x \geq M} x^i f(x) d(x)$ of order i smaller or equal to 2, he obtains the coefficients $A_1 = \exp\left\{2\left[\mu_1(Z) + 0.6931 B_1\right]\right\}$, $B_1 = 1.7099\left\{0.5\mu_2(Z) - \left[\mu_1(Z)\right]^2\right\}^{0.5}$, $A_2^2 = \dfrac{M_2(Y) - 2\left[M_1(Y)\right]^2}{0.6840}$, and $B_2 = 2\left[M_1(Y) - 0.6931 A_2\right]$ where $Z = \begin{cases} \ln(Y) & \text{pour } Y < M \\ 0 & \text{pour } Y \geq M \end{cases}$, $Y$ the subjacent distribution, $\mu_i(Z)$ the ith moment of $Z$, and $M_i(Y)$ the ith partial moment of $Y$. Shore's method leads to better estimates of the percentiles than those obtained by Clements's method [4], insofar as the expected values being similar, the mean squared error is much lower by Shore's method [19,20,21].

## 5. Distribution and moments of the estimators of the indices suggested for a normal distribution.

The studied characteristic of the process is supposed normally distributed with mean $\mu$ and variance $\sigma^2$. Two natural estimators of $C_{pu}(u,v)$ can be considered, differing in the way the variance $\sigma^2$ is estimated. We define the estimators $\hat{C}_{pu,n}(u,v)$ and $\hat{C}_{pu,n-1}(u,v)$ as

$$\hat{C}_{pu,n}(u,v) = \frac{U - T - u\hat{A}_u^*}{3\sqrt{S_n^2 + v\hat{A}_u^{*2}}}, \text{ and } \hat{C}_{pu,n-1}(u,v) = \frac{U - T - u\hat{A}_u^*}{3\sqrt{S_{n-1}^2 + v\hat{A}_u^{*2}}},$$

where $\hat{A}_u^* = \max\left\{\overline{X} - T, (T - \overline{X})/k\right\}$, $\overline{X} = \left(\sum_{i=1}^n X_i\right)/n$, $S_n^2 = \sum_{i=1}^n (X_i - \overline{X})^2/n$, and

$S_{n-1}^2 = \sum_{i=1}^n (X_i - \overline{X})^2/(n-1)$.

The two estimators are related by $\hat{C}_{pu,n-1}(u,v) = ((n-1)/n)^{1/2} \hat{C}_{pu,n}(u, v(n-1)/n)$. The study of the statistical properties of these estimators is facilitated insofar as $C_{pu}(u,v)$ can be expressed according to Chen's and Pearn's $C_p^*(u,v)$ index [4]. Let us recall that for an interval $[L;U]$ and a target $T$ not centered, they introduce the family

$$C_p^*(u,v) = \frac{d^* - uA^*}{3\sqrt{\sigma^2 + vA^2}}, \; u \geq 0, \; v \geq 0,$$

where $A^* = \max\left\{d^*(\mu - T)/D_u, d^*(T - \mu)/D_l\right\}$,

$A = \max\left\{d(\mu - T)/D_u, d(T - \mu)/D_l\right\} = dA^*/d^*$,

$D_u = U - T$, $D_l = T - L$, and $d^* = \min\{D_u, D_l\}$.

In the case of a single upper tolerance $U$, considering that the risk of a deviation towards the left of the target is $k$ time less serious than towards the upper limit $U$ can be interpreted as the positioning of a lower limit $L$ so that $T - L = D_l = kD_u$. In these conditions $d = (1+k)D_u/2$, $d^* = D_u$, $A_u^* = A^* = \max\left\{\mu - T, (T - \mu)/k\right\}$, $A = dA^*/d^* = (1+k)A^*/2$, and

$$C_{pu}(u,v) = \frac{D_u - u A_u^*}{3\sqrt{\sigma^2 + v A_u^{*2}}} = \frac{d^* - uA^*}{3\sqrt{\sigma^2 + v(d^*/d)^2 A^2}} = C_p''\left(u, v(d^*/d)^2\right)$$

$$= C_p''\left(u, 4v/(1+k)^2\right). \tag{1}$$

Assuming that $\delta = \sqrt{n}(\mu - T)/\sigma$, $\lambda = \delta^2$, and $D^* = \sqrt{n}d^*/\sigma$, Grau [8] gives the r-th moment of $\hat{C}_{p,n}^*(u,v)$ in the form

$$E\left(\hat{C}_{p,n}^*(u,v)\right)^r = \frac{1}{3^r}\sum_{i=0}^r (-u)^i \binom{r}{i}\left(\frac{D^*}{\sqrt{2}}\right)^{r-i} \frac{e^{-\lambda/2}}{2\sqrt{\pi}} \times \sum_{j=0}^\infty \frac{\delta^j 2^{j/2}}{j!} \frac{\Gamma(c-a)\Gamma(b)}{\Gamma(c)} \gamma_n(i,j), \tag{2}$$

with $\gamma_n(i,j) = \left[\left(\frac{d^*}{D_u}\right)^i {}_2F_1(a,b;c;z_u) + (-1)^j \left(\frac{d^*}{D_l}\right)^i {}_2F_1(a,b;c;z_l)\right]$,

where $(-u)^i$ should be interpreted as 1 when $i = 0$, also for the case $u = 0$, and ${}_2F_1(a,b;c;z)$ is the Gaussian hypergeometric function (Abramowitz and Stegun [1]) with parameters $a = r/2$, $b = (1+i+j)/2$, $c = (n+i+j)/2$, $z_u = 1 - (d/D_u)^2 v$, and $z_l = 1 - (d/D_l)^2 v$.

In order to distinguish the properties of the estimators of the indices $C_{pu}(u,v)$ and $C_{pl}(u,v)$ subsequently, we assume $B_u = \sqrt{n}(U - T)/\sigma$, which in that case is equal to $D^*$. From relations (1) and (2) we deduce the r-th moment of $\hat{C}_{pu,n}(u,v)$,

$$E\left(\hat{C}_{pu,n}(u,v)\right)^r = \frac{1}{3^r}\sum_{i=0}^r (-u)^i \binom{r}{i}\left(\frac{B_u}{\sqrt{2}}\right)^{r-i} \frac{e^{-\lambda/2}}{2\sqrt{\pi}} \times \sum_{j=0}^\infty \frac{\delta^j 2^{j/2}}{j!} \frac{\Gamma(c-a)\Gamma(b)}{\Gamma(c)} \gamma_{u,n}(i,j), \tag{3}$$

with $\gamma_{u,n}(i,j) = \left[ {}_2F_1\left(a,b;c;1-v\right) + (-1)^j k^{-i} {}_2F_1\left(a,b;c;1-v/k^2\right)\right]$.

In the same way we obtain $E\left(\hat{C}_{pl,n}(u,v)\right)^r$ by replacing $B_u$ by $B_l = \sqrt{n}\left(T-L\right)/\sigma$ and $\gamma_{u,n}(i,j)$ by $\gamma_{l,n}(i,j) = \left[ k^{-i} {}_2F_1\left(a,b;c;1-v/k^2\right) + (-1)^j {}_2F_1\left(a,b;c;1-v\right)\right]$.

## 5.1 Estimation and distribution of $C_{pu}$ and $C_{pl}$

As previously, if we consider that the choice of the constant $k$ can be interpreted as the positioning of a lower limit $L$ so that $T-L = D_l = kD_u$, then $C_{pu}$ can be expressed according to the usual index $C_p$ by the relation $C_{pu} = \left(2/(1+k)\right)C_p$. From the density of probability and moments of $\hat{C}_{p,n-1}$ given by Kotz and Lovelace [10], we obtain

$$f_{\hat{C}_{pu,n-1}}(x) = \frac{(n-1)^{(n-1)/2}}{C_{pu}\Gamma\left((n-1)/2\right)2^{(n-3)/2}}\left(\frac{C_{pu}}{x}\right)^n e^{-((n-1)/2)\left(C_{pu}/x\right)^2}, \quad x > 0,$$

and $E\left(\hat{C}_{pu,n-1}\right)^r = \left((n-1)/2\right)^{r/2}\dfrac{\Gamma\left((n-1-r)/2\right)}{\Gamma\left((n-1)/2\right)}C_{pu}^r$.

In particular $E\left(\hat{C}_{pu,n-1}\right) = b_f^{-1}C_{pu}$ and $V\left(\hat{C}_{pu,n-1}\right) = \left(\dfrac{n-1}{n-3} - b_f^{-2}\right)C_{pu}^2$ where $b_f = \sqrt{\dfrac{2}{n-1}}\dfrac{\Gamma\left((n-1)/2\right)}{\Gamma\left((n-2)/2\right)}$.

Moreover, since $\hat{C}_{pu,n} = \left(n/(n-1)\right)^{1/2}\hat{C}_{pu,n-1}$,

$$f_{\hat{C}_{pu,n}}(x) = \frac{n^{(n-1)/2}}{C_{pu}\Gamma\left((n-1)/2\right)2^{(n-3)/2}}\left(\frac{C_{pu}}{x}\right)^n e^{-(n/2)\left(C_{pu}/x\right)^2}, \quad x > 0,$$

and $E\left(\hat{C}_{pu,n}\right)^r = \left(n/2\right)^{r/2}\dfrac{\Gamma\left((n-1-r)/2\right)}{\Gamma\left((n-1)/2\right)}C_{pu}^r$.

In particular $E\left(\hat{C}_{pu,n}\right) = c_f^{-1}C_{pu}$ and $V\left(\hat{C}_{pu,n}\right) = \left(\dfrac{n}{n-3} - c_f^{-2}\right)C_{pu}^2$ where $c_f = \sqrt{\dfrac{2}{n}}\dfrac{\Gamma\left((n-1)/2\right)}{\Gamma\left((n-2)/2\right)}$.

In a similar way we obtain the density and the moments of $\hat{C}_{pl,n}$ and $\hat{C}_{pl,n-1}$ replacing $C_{pu}$ by $C_{pl}$.

## 5.2 Estimation and distribution of $C_{pku}$ and $C_{pkl}$

If we consider that the choice of the constant $k$ can be interpreted as the positioning of a lower limit $L$ so that $T-L = D_l = kD_u$, then $C_{pku} = C_{pu}(1,0) = C_p^{"}(1,0) = C_{pk}^{"}$ according to (1). Since $D^* = B_u$, assuming that $D = n^{1/2}d/\sigma$, according to the appendix, we obtain

$$f_{\hat{C}_{pku,n}}(x) = \begin{cases} -\int_1^\infty J_1(x,t)dt, & x < 0 \\ \int_0^1 J_1(x,t)dt, & x > 0 \end{cases},$$

where $J_1(x,t) = f_K\left(\left(\dfrac{B_u\left(1-\sqrt{t}\right)}{3x}\right)^2\right) f_{Y_u}\left(B_u^2 t\right) \dfrac{2B_u^2}{x}\left(\dfrac{B_u\left(1-\sqrt{t}\right)}{3x}\right)^2$, and

$f_{Y_u}(y) = \dfrac{1}{2\sqrt{y}}\left(\phi\left(\sqrt{y}-\delta\right) + k\phi\left(k\sqrt{y}+\delta\right)\right)$ when y > 0, with $\phi(x)$ the probability density of a N(0,1) distribution. By substitution, the density can still be written

$$f_{\hat{C}_{pku,n}}(x) = \begin{cases} -\int_{B_u^2}^{\infty} J_1^{\cdot}(x,t)dt, & x < 0 \\ \int_0^{B_u^2} J_1^{\cdot}(x,t)dt, & x > 0 \end{cases},$$

where $J_1^{\cdot}(x,t) = f_K\left(\left(\dfrac{B_u-\sqrt{t}}{3x}\right)^2\right) f_{Y_u}(t)\dfrac{2}{x}\left(\dfrac{B_u-\sqrt{t}}{3x}\right)^2$.

Since $\hat{C}_{pku,n-1} = \left((n-1)/n\right)^{1/2}\hat{C}_{pku,n}$, we deduce

$$f_{\hat{C}_{pku,n-1}}(x) = \begin{cases} -\int_1^{\infty} L_1(x,t)dt, & x < 0 \\ \int_0^1 L_1(x,t)dt, & x > 0 \end{cases},$$

where $L_1(x,t) = f_K\left(\dfrac{n-1}{n}\left(\dfrac{B_u\left(1-\sqrt{t}\right)}{3x}\right)^2\right) f_{Y_u}\left(B_u^2 t\right)\dfrac{2B_u^2}{x}\dfrac{n-1}{n}\left(\dfrac{B_u\left(1-\sqrt{t}\right)}{3x}\right)^2$.

Since $D^* = B_u$, we find the density of $\hat{C}_{pk,n-1}^* = \hat{C}_{pku,n-1}$ given by Pearn, Lin and Chen [18]. Only the moments of order 1 and 2 of $\hat{C}_{pk,n-1}^*$ are explicitly given by Pearn and Chen [13]. Thus, we use the equation (3) for the r-th moment which, in addition, leads to expressions of the first two moments simpler than those given by Pearn and Chen [13].

$$E\left(\hat{C}_{pku,n}\right)^r = \dfrac{1}{3^r}\sum_{i=0}^{r}(-1)^i\binom{r}{i}\left(\dfrac{B_u}{\sqrt{2}}\right)^{r-i}\dfrac{e^{-\lambda/2}}{2\sqrt{\pi}}\times\dfrac{\Gamma\left((n-r-1)/2\right)}{\Gamma\left((n-1)/2\right)}$$

$$\times\sum_{j=0}^{\infty}\dfrac{\delta^j 2^{j/2}}{j!}\Gamma\left((1+i+j)/2\right)\gamma_u(i,j),$$

where $\gamma_u(i,j) = 1 + (-1)^j k^{-i}$. In particular, Grau [8],

$$E\left(\hat{C}_{pku,n}\right) = c_f^{-1}\left[C_{pku} + \dfrac{1+k^{-1}}{3\sqrt{n}}\left(|\delta|\Phi\left(-|\delta|\right) - \dfrac{e^{-\delta^2/2}}{\sqrt{2\pi}}\right)\right],$$

$$E\left(\hat{C}_{pku,n}\right)^2 = \dfrac{n}{(n-3)}\left[C_{pku}^2 - \dfrac{2(U-T)(1+k^{-1})}{9\sigma\sqrt{n}}\left(\dfrac{e^{-\delta^2/2}}{\sqrt{2\pi}} - |\delta|\Phi\left(-|\delta|\right)\right)\right.$$

$$\left. + \dfrac{1+k^{-2}}{18n} + \dfrac{\delta\left(1-k^{-2}\right)}{9n}\left(\dfrac{e^{-\delta^2/2}}{\sqrt{2\pi}} - |\delta|\Phi\left(-|\delta|\right) + \dfrac{1}{2\delta}\left\{1-2\Phi\left(-\delta\right)\right\}\right)\right].$$

$E\left(\hat{C}_{pku,n-1}\right)^r$, $E\left(\hat{C}_{pku,n-1}\right)$ and $E\left(\hat{C}_{pku,n-1}\right)^2$ are obtained without difficulty since

$\left(\hat{C}_{pku,n-1}\right)^r = \left((n-1)/n\right)^{r/2}\left(\hat{C}_{pku,n}\right)^r$.

For $\hat{C}_{pkl,n}$ and $\hat{C}_{pkl,n-1}$ we obtain similar results replacing $f_{Yu}(y)$ by

$$f_{Yl}(y) = \frac{1}{2\sqrt{y}}\left(k\phi\left(k\sqrt{y}-\delta\right)+\phi\left(\sqrt{y}+\delta\right)\right), \ B_u \text{ by } B_l, \text{ and } \gamma_u(i,j) \text{ by } \gamma_l(i,j) = k^{-i}+(-1)^j.$$

## 5.3 Estimation and distribution of $C_{pmu}$ and $C_{pml}$

If we consider that the choice of the constant $k$ can be interpreted as the positioning of a lower limit $L$ so that $T-L=D_l=kD_u$, then $C_{pmu}=C_{pu}(0,1)=C_p''\left(0,\left(d^*/d\right)^2\right)$ according to (1).

From the appendix,

$$f_{\hat{C}_{pmu,n}}(x) = \int_0^1 J_2(x,t)dt, \qquad x > 0,$$

where $J_2(x,t)=f_K\left(G(x)(1-t)\right)f_{Yu}\left(G(x)t\right)2x^{-1}\left(G(x)\right)^2$, and $G(x)=\left(B_u/(3x)\right)^2$, or

$$f_{\hat{C}_{pmu,n}}(x) = \int_0^{G(x)} J_2'(x,t)dt, \qquad x > 0,$$

where $J_2'(x,t)=f_K\left(G(x)-t\right)f_{Yu}(t)2x^{-1}G(x)$.

Since $\hat{C}_{pmu,n-1}=\left((n-1)/n\right)^{1/2}\hat{C}_{pu,n}\left(0,(n-1)/n\right)$, we deduce

$$f_{\hat{C}_{pmu,n-1}}(x) = \int_0^1 L_2(x,t)dt, \qquad x > 0$$

where $L_2(x,t)=f_K\left((n-1/n)G(x)(1-t)\right)f_{Yu}\left(G(x)t\right)2(n-1/n)x^{-1}\left(G(x)\right)^2$.

The moments are obtained from the equation (3):

$$E\left(\hat{C}_{pmu,n}\right)^r = \frac{1}{3^r}\left(\frac{B_u}{\sqrt{2}}\right)^r\frac{e^{-\lambda/2}}{2\sqrt{\pi}}\times\sum_{j=0}^\infty\frac{\delta^j 2^{j/2}}{j!}\frac{\Gamma(c-a)\Gamma(b)}{\Gamma(c)}\gamma_{u,n}(j),$$

$$E\left(\hat{C}_{pmu,n-1}\right)^r = \left(\frac{n-1}{n}\right)^{r/2}\frac{1}{3^r}\left(\frac{B_u}{\sqrt{2}}\right)^r\frac{e^{-\lambda/2}}{2\sqrt{\pi}}\times\sum_{j=0}^\infty\frac{\delta^j 2^{j/2}}{j!}\frac{\Gamma(c-a)\Gamma(b)}{\Gamma(c)}\gamma_{u,n-1}(j),$$

where $a=r/2$, $b=(1+j)/2$, $c=(n+j)/2$, $\gamma_{u,n}(j)=1+(-1)^j \ _2F_1\left(a,b;c;1-k^{-2}\right)$, and

$\gamma_{u,n-1}(j)=\ _2F_1\left(a,b;c;n^{-1}\right)+(-1)^j \ _2F_1\left(a,b;c;1-(n-1)n^{-1}k^{-2}\right)$.

For $\hat{C}_{pml,n}$ and $\hat{C}_{pml,n-1}$, we obtain similar results replacing $f_{Yu}(y)$ by $f_{Yl}(y)$, $B_u$ by $B_l$,

$\gamma_{u,n}(j)$ by $\gamma_{l,n}(j)=\ _2F_1\left(a,b;c;1-k^{-2}\right)+(-1)^j$, and $\gamma_{u,n-1}(j)$ by

$\gamma_{l,n-1}(j)=\ _2F_1\left(a,b;c;1-(n-1)n^{-1}k^{-2}\right)+(-1)^j \ _2F_1\left(a,b;c;n^{-1}\right)$.

## 5.4 Estimation and distribution of $C_{pmku}$ and $C_{pmkl}$

If we consider that the choice of the constant $k$ can be interpreted as the positioning of a lower limit $L$ so that $T-L=D_l=kD_u$, then $C_{pmku}=C_{pu}(1,1)=C_p''\left(1,\left(d^*/d\right)^2\right)$. From the

appendix, we obtain

$$f_{\hat{C}_{pmku,n}}(x) = \begin{cases} -\int_1^\infty J(x,t)dt, & -\frac{1}{3} < x < 0 \\ \int_0^1 J(x,t)dt, & x > 0 \end{cases}$$

where $J(x,t) = f_K\left(\left(\dfrac{B_u - \sqrt{tH(x)}}{3x}\right)^2 - tH(x)\right) f_{Yu}(tH(x)) \dfrac{2H(x)}{x}\left(\dfrac{B_u - \sqrt{tH(x)}}{3x}\right)^2$, and

$H(x) = \left[B_u \big/ (1+3x)\right]^2$, or

$$f_{\hat{C}_{pmku,n}}(x) = \begin{cases} -\displaystyle\int_{H(x)}^{\infty} J'(x,t)\,dt, & -\dfrac{1}{3} < x < 0 \\[2mm] \displaystyle\int_0^{H(x)} J'(x,t)\,dt, & x > 0 \end{cases}$$

where $J'(x,t) = f_K\left(\left(\dfrac{B_u - \sqrt{y}}{3x}\right)^2 - y\right) f_{Yu}(y) \dfrac{2}{x}\left(\dfrac{B_u - \sqrt{y}}{3x}\right)^2$.

Since $\hat{C}_{pmku,n-1} = \left((n-1)/n\right)^{1/2} \hat{C}_{pu,n}\left(1,(n-1)/n\right)$, we deduce

$$f_{\hat{C}_{pmku,n-1}}(x) = \begin{cases} -\displaystyle\int_1^{\infty} L(x,t)\,dt, & -\dfrac{1}{3}\left(\dfrac{n}{n-1}\right)^{1/2} < x < 0 \\[2mm] \displaystyle\int_0^1 L(x,t)\,dt, & x > 0 \end{cases}$$

where $L(x,t) = f_K\left(\dfrac{n-1}{n}\left(\left(\dfrac{B_u - \sqrt{tH(x)}}{3x}\right)^2 - tH(x)\right)\right) f_{Yu}(tH(x)) \dfrac{2H(x)}{x}\dfrac{n-1}{n}\left(\dfrac{B_u - \sqrt{tH(x)}}{3x}\right)^2$.

The moments are obtained from the equation (3):

$$E\left(\hat{C}_{pmku,n}\right)^r = \dfrac{1}{3^r}\sum_{i=0}^{r}(-1)^i\binom{r}{i}\left(\dfrac{B_u}{\sqrt{2}}\right)^{r-i}\dfrac{e^{-\lambda/2}}{2\sqrt{\pi}}\times\sum_{j=0}^{\infty}\dfrac{\delta^j 2^{j/2}}{j!}\dfrac{\Gamma(c-a)\Gamma(b)}{\Gamma(c)}\gamma_{u,n}(i,j),$$

$$E\left(\hat{C}_{pmku,n-1}\right)^r = \left(\dfrac{n-1}{n}\right)^{r/2}\dfrac{1}{3^r}\sum_{i=0}^{r}(-1)^i\binom{r}{i}\left(\dfrac{B_u}{\sqrt{2}}\right)^{r-i}\dfrac{e^{-\lambda/2}}{2\sqrt{\pi}}\times\sum_{j=0}^{\infty}\dfrac{\delta^j 2^{j/2}}{j!}\dfrac{\Gamma(c-a)\Gamma(b)}{\Gamma(c)}\gamma_{u,n-1}(i,j),$$

where $a = r/2$, $b = (1+i+j)/2$, $c = (n+i+j)/2$, $\gamma_{u,n}(i,j) = \left[1+(-1)^j k^{-i}\,{}_2F_1\left(a,b;c;1-k^{-2}\right)\right]$,

and $\gamma_{u,n-1}(i,j) = {}_2F_1\left(a,b;c;n^{-1}\right) + (-1)^j k^{-i}\,{}_2F_1\left(a,b;c;1-(n-1)n^{-1}k^{-2}\right)$.

For $\hat{C}_{pmkl,n}$ and $\hat{C}_{pmkl,n-1}$, we obtain similar results replacing $f_{Yu}(y)$ by $f_{Yl}(y)$, $B_u$ by $B_l$, $\gamma_{u,n}(j)$ by $\gamma_{l,n}(i,j) = \left[k^{-i}\,{}_2F_1\left(a,b;c;1-k^{-2}\right)+(-1)^j\right]$, and $\gamma_{u,n-1}(i,j)$ by $\gamma_{l,n-1}(i,j) = k^{-i}\,{}_2F_1\left(a,b;c;1-(n-1)n^{-1}k^{-2}\right)+(-1)^j\,{}_2F_1\left(a,b;c;n^{-1}\right)$.

## 6  Example

To illustrate how these indices can be used, we present a study carried out within the chemistry company Atofina (France). This company manufactures polymer granulates. Before polymerization, an additive, the maleic anhydride, is blended. The target for this additive is 480 ppm and the necessary minimal value is 400 ppm, which represents the lower limit of tolerance. The cost of this additive not being very important, no upper limit has been fixed. Analyses are made regularly on the granulates at the end of the process. 86 values of this variable have been collected and are given below :

730 730 740 820 750 630 670 630 750 780 790 890 670 675 710 720 920 620 620
710 680 750 660 780 830 850 645 645 710 815 695 645 710 815 695 840 790 820

660 700 685 490 610 685 590 620 700 730 750 780 730 770 760 780 640 830 785
780 680 650 585 650 600 610 640 640 575 630 610 665 690 650 640 645 630 700
800 685 650 675 640 655 610 625 580 610

Figure 3 gives the histogram of the 86 observations studied. A khi square test carried out on the classes of this histogram gives a p-value of 0.003. This p-value leads us to reject the normality of the observations. Thus we use the indices defined for a non-normal distribution and estimate the quantiles by the method of Shore [18]. We obtain $\hat{M} = 685$, $A_1 = 677.715$, $B_1 = 0.050$, $A_2 = 48.307$, $B_2 = 695.712$, from which $\hat{L}_{0.00135} = 486.606$, $\hat{U}_{0.99865} = 1014.840$, and $\hat{C}_{pl} = 0.40$. The low value of $\hat{C}_{pl}$ indicates that the process is potentially not capable because the dispersion is too important according to the interval $[L;T]$. The histogram shows however that no unit is below the tolerance, but that the median estimated at 685 ppm is largely above the target $T = 480$ ppm. However, even if the company assesses that this deviation is less serious than in the other direction, the resulting additional cost requires the use of indices of real capability taking into account the position of the median according to the target. As an example we will consider the 2 situations where the gravity was judged 3 times ($k = 3$), then 10 times ($k = 10$) less important. For $k=3$, we obtain $\hat{C}_{pkl} = 0.06$, $\hat{C}_{pml} = 0.28$, and $\hat{C}_{pmkl} = 0.04$, and for $k = 10$, $\hat{C}_{pkl} = 0.30$, $\hat{C}_{pml} = 0.39$, and $\hat{C}_{pmkl} = 0.29$. Of course, when the deviation takes place on the right, the smaller the risk is, the closer to the potential capability the indices of real capability are.

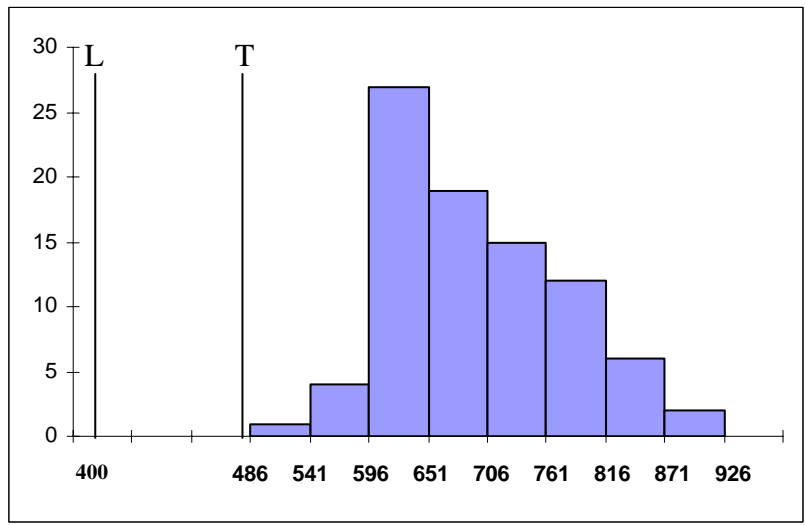

Figure 3. Histogram of the data

## 7 Conclusion

In this paper we first review the existing generalizations of the basic capability indices $C_p$, $C_{pk}$, $C_{pm}$, and $C_{pmk}$ which have been proposed to handle one-sided tolerances processes. In this situation, the risk of a deviation towards the side opposed to the limit of tolerance is not considered very serious and thus, in general, not quantified. However, for problems of cost among other things, this risk cannot be completely ignored. When a ratio is established, at least approximately, between the risks taken for deviations to the left and to the right of the target, we introduce indices of capability preserving the properties and interpretations of the usual family $C_p$, $C_{pk}$, $C_{pm}$, and $C_{pmk}$. For a process of normal distribution, the densities and moments of the estimators are given. They ensue from the properties established for the family $C_p^{"}(u,v)$, used in the case of asymmetrical tolerances. When the distribution of the process is not normal, we propose simple extensions using the median and the quantiles of order 0.00135 and 0.99865. Finally an example shows the way to calculate these quantities on a real example of non-normal data, using the method of Shore [18] in order to estimate the extreme quantiles.

## Appendix

Let us consider the case of asymmetrical tolerances where $D_u = U - T$, $D_l = T - L$, and $d^* = \min\{D_u, D_l\}$. Assuming that $D^* = n^{1/2} d^*/\sigma$, $u' = d^*/d$, $\delta = \sqrt{n}(\mu - T)/\sigma$, $S(x) = \left[D^*/(u'+3x)\right]^2$, $Z = n^{1/2}(\bar{X} - T)/\sigma$, $K = nS_n^2/\sigma^2$ the probability density of which $f_K(x)$ is a $\chi_{n-1}^2$, $\phi(x)$ the probability density of a N(0,1), and $Y = \left[\max\{(d/D_u)Z, -(d/D_l)Z\}\right]^2$ the probability density of which is

$$f_Y(y) = \frac{1}{2\sqrt{y}}\left(\frac{D_u}{d}\phi\left(\frac{D_u}{d}\sqrt{y} - \delta\right) + \frac{D_l}{d}\phi\left(\frac{D_l}{d}\sqrt{y} + \delta\right)\right) \text{ when y} > 0.$$

For $\hat{A}^* = \max\{d^*(\bar{X} - T)/D_u, d^*(T - \bar{X})/D_l\}$ and $\hat{A} = \max\{d(\bar{X} - T)/D_u, d(T - \bar{X})/D_l\} = d\hat{A}^*/d^*$, Pearn, Lin and Chen [17] obtain the distribution of the estimator

$$\hat{C}_{pmk,n}^{"} = \frac{d^* - \hat{A}^*}{3\sqrt{S_n^2 + \hat{A}^2}} = \frac{D^* - u'\sqrt{Y}}{3\sqrt{K + Y}} \text{ of } C_{pmk}^{"} \text{ in the form}$$

$$f_{\hat{C}_{pmk,n}^{"}}(x) = \begin{cases} -\int_1^\infty I(x,t)dt, & -\frac{u'}{3} < x < 0 \\ \int_0^1 I(x,t)dt, & x > 0 \end{cases},$$

where $I(x,t) = f_K\left(\left(\frac{D^* - u'\sqrt{tS(x)}}{3x}\right)^2 - tS(x)\right)f_Y\left(tS(x)\right)\frac{2S(x)}{x}\left(\frac{D^* - u'\sqrt{tS(x)}}{3x}\right)^2$.

More generally, let us consider the index $C_p^{"}(u,v) = \frac{d^* - uA^*}{3\sqrt{\sigma^2 + vA^2}}$ when $u, v > 0$. Using the same demonstration as Pearn, Lin and Chen [17], and assuming that $h = u\,d^*/d$, and $K(x) = \left[D^*/\left(h + 3x\sqrt{v}\right)\right]^2$, we obtain the density of probability of the estimator $\hat{C}_{p,n}^{"}(u,v)$ of $C_p^{"}(u,v)$ in the form

$$f_{\hat{C}_{p,n}^{"}(u,v)}(x) = \begin{cases} -\int_1^\infty J(x,t)dt, & -\frac{h}{3\sqrt{v}} < x < 0 \\ \int_0^1 J(x,t)dt, & x > 0 \end{cases},$$

where $J(x,t) = f_K\left(\left(\frac{D^* - h\sqrt{tK(x)}}{3x}\right)^2 - vtK(x)\right)f_Y\left(tK(x)\right)\frac{2K(x)}{x}\left(\frac{D^* - h\sqrt{tK(x)}}{3x}\right)^2$.

Let us note that when $u > 0$ and $v = 0$, the previous expression is still valid if we consider that $-h/\left(3\sqrt{v}\right) = -\infty$. For $u = 0$ and $v > 0$, the expression is also valid when x > 0, the density being null when x < 0.

In the particular case where $u > 0$, $v = 0$, we obtain $K(x) = \left[D^*/h\right]^2 = \left[D/u\right]^2$, with $D = n^{1/2}\,d/\sigma$, and

$$f_{\hat{C}_{p,n}^{'}(u,0)}(x) = \begin{cases} -\int_1^\infty J_1(x,t)dt, & x < 0 \\ \int_0^1 J_1(x,t)dt, & x > 0 \end{cases},$$

where $J_1(x,t) = f_K\left(\left(\dfrac{D^*\left(1-\sqrt{t}\right)}{3x}\right)^2\right)f_Y\left(\left(\dfrac{D}{u}\right)^2 t\right)\dfrac{2}{x}\left(\dfrac{D}{u}\right)^2\left(\dfrac{D^*\left(1-\sqrt{t}\right)}{3x}\right)^2.$

In the particular case where $u = 0$, $v > 0$, we obtain $K(x) = \left[D^*/\left(3x\sqrt{v}\right)\right]^2$, and

$$f_{\hat{C}_{p,n}^{'}(0,v)}(x) = \int_0^1 J_2(x,t)dt, \qquad x > 0,$$

where $J_2(x,t) = f_K\left(\left(\dfrac{D^*}{3x}\right)^2\left(1-t\right)\right)f_Y\left(\left(\dfrac{D^*}{3x}\right)^2\dfrac{t}{v}\right)\dfrac{2}{xv}\left(\dfrac{D^*}{3x}\right)^4.$

**References**

[1]    Abramowitz, M., and Stegun, I. A. (1965). Handbook of Mathematical functions with formulas, graphs, and mathematical tables. Dover publications, New York.

[2]    Chan, L.K., Cheng, S.W. and Spiring, F.A. (1988). A new measure of process capability : C$_{pm}$. *Journal of Quality Technology*, *20*, 162–175.

[3]    Chen, K. S., Pearn, W. L. and Lin, P. C. (1999). A new generalization of $C_{pm}$ for processes with asymmetric tolerances. *International Journal of Reliability, Quality and Safety Engineering*, 6(4), 383-398.

[4]    Chen, K. S. and Pearn, W. L. (2001). Capability indices for processes with asymmetric tolerances. *Journal of the Chinese Institute of Engineers*, 24(5), 559-568.

[5]    Clements, J. A. (1989). Process capability calculations for non-normal distributions. *Quality Progress*, 95-100.

[6]    Ding, J. (2004). A method of estimating the process capability index from the first four moments of non-normal data. *Quality and Reliability Engineering International*, 20, 787-805.

[7]    Grau, D. (2005). New capability indices for non-normal processes. *Submitted*.

[8]    Grau, D. (2006). On the choice of a capability index for asymmetric tolerances. *Submitted*.

[9]    Kane, V. E. (1986). Process Capability Indices. *Journal of Quality Technology*, 18, 41-52.

[10]    Kotz, S. and Lovelace, C.R. (1998). *Process Capability Indices in Theory and Practice*, Arnold, London.

[11]    McCormack, D.W., Harris, I.R., Hurwitz, A.M. and Spagon, P.D. (2000). Capability indices for non-normal data. *Quality Engineering*, 12, 489-495.

[12]    Pearn, W. L. and Chen, K. S. (1995). Estimating process capability indices for non-Normal Pearsonian populations. *Quality and Reliability Engineering International*, 11(15), 386-388.

[13]   Pearn, W. L. and Chen, K. S. (1998). New generalization of process capability index $C_{pk}$. *Journal of Applied Statistics*, 25(6), 801-810.

[14]   Pearn, W. L., Chen, K. S. and Lin, G. H. (1999). A generalization of Clements' method for non-Normal Pearsonian processes with asymmetric tolerances. *International Journal of Quality and Reliability Management*, 16(5), 507-521.

[15]   Pearn, W. L., Kotz, S.(1994). Application of Clements' method for calculating second and third generation process capability indices for Non-normal personian populations. *Quality Engineering*, 7, 139-145.

[16]   Pearn, W. L., Lin, P. C. and Chen, K. S. (1999). On the generalizations of the capability index $C_{pmk}$ for asymmetric tolerances, *Far East Journal of Theoretical Statistics*, 3, 49-66.

[17]   Pearn, W. L., Lin, P.C. and Chen, K. S. (2001). Estimating process capability index $C_{pmk}^{"}$ for asymmetric tolerances: distributional properties, *Metrika*, 54, 261-279.

[18]   Pearn, W. L., Lin, P.C. and Chen, K. S. (2004). The $C_{pk}^{"}$ index for asymmetric tolerances: implications and inference, *Metrika*, 60, 119-136.

[19]   Shore, H. (1998). A new approach to analysing non-normal quality data with application to process capability analysis. *International Journal of Production Research*, 36(7), 1917-1933.

[20]   Shore, H. (1998). Approximating an unknown distribution when distribution information is extremely limited. *Communications in statistics (Simulation and computation)*, 27(2), 501-523.

[21]   Shore, H. (2004). Non-normal populations in quality applications: a revisited perspective. *Quality and Reliability Engineering International*, 20, 375-382.

[22]   Tang, L.C., and Than, S.E.E (1999). Computing process capability indices for non-normal data: A review and comparative study. *Quality and Reliability Engineering International*, 15, 339-353.

[23]   Vännman, K. (1998). Families of capability indices for one-sided specification limits. *Statistics*, 31, 43-66.